\documentclass[11pt]{amsart}
\usepackage{amsmath}
\usepackage{amsfonts}

\usepackage{amscd}
\usepackage{amsthm}
\usepackage{amssymb} \usepackage{latexsym}
\usepackage{eufrak}
\usepackage{euscript}
\usepackage{epsfig}
\usepackage{graphics}
\usepackage{array}
\usepackage{enumerate}
\usepackage{dsfont}
\usepackage{color}
\usepackage{wasysym}
\usepackage{hyperref}
\usepackage{pdfsync}

\newcommand{\bel}[1]{\begin{equation}\label{#1}}

\newcommand{\be}{\begin{equation}}

\newcommand{\ba}{\begin{eqnarray}}
\newcommand{\ea}{\end{eqnarray}}

\newcommand{\qe}{\end{equation}}

\newcommand{\Hmm}[1]{\leavevmode{\marginpar{\tiny%
$\hbox to 0mm{\hspace*{-0.5mm}$\leftarrow$\hss}%
\vcenter{\vrule depth 0.1mm height 0.1mm width \the\marginparwidth}%
\hbox to
0mm{\hss$\rightarrow$\hspace*{-0.5mm}}$\\\relax\raggedright #1}}}

\theoremstyle{theorem}
\newtheorem{thm}{Theorem}[section]

\theoremstyle{example}

\theoremstyle{corollary}

\theoremstyle{lemma}
\newtheorem{lem}[thm]{Lemma}
\theoremstyle{definition}
\newtheorem{defi}[thm]{Definition}
\theoremstyle{proof}
\theoremstyle{remark}
\newtheorem{rem}[thm]{Remark}

\begin{document}

\title[Properties for CD Inequalities with Unbounded Laplacians]{Equivalent Properties for CD Inequalities on Graphs with Unbounded Laplacians}
\author{Chao Gong}
\email{elfmetal@ruc.eud.cn}
\address{Department of Mathematics,Information School,
Renmin University of China,
Beijing 100872, China
}

\author{Yong Lin}
\email{linyong01@ruc.edu.cn}
\address{Department of Mathematics,Information School,
Renmin University of China,
Beijing 100872, China}

\begin{abstract}The CD equalities were introduced to imply the gradient estimate of laplace operator on graphs. This article is based on the unbounded Laplacians, and finally concludes some equivalent properties of the CD(K,$\infty$)and CD(K,n).
\end{abstract}
\maketitle

\section{Introduction}

Graph theory is the basic theory of the study of graphs and networks. The spectral graph theory, which is used for describing the structure and characteristic of graphs by adjacency matrix or the spectral density of Laplacian matrix, is the classic method for studying graphs described in {\cite{9}}.

\indent We already know that we can find the curvature by solving a partial differential equation, and there are more examples for the geometric analysis such as the famous Li-Yau gradient estimate. Moreover, we can use some data to describe the graphs and optimize it such as the Cheeger constant on graphs.\\
\indent The Laplacians on graph have always been an important research topic. In fact, Laplacians can be seen as the discrete analog form of the Schrodinger operator, or as the generator of symmetric Markov process. Laplacians always appear in the topics on the research of discrete structure for heat equations like in {\cite{2}}.\\
\indent As for the Laplacians on graphs, the properties are different on different occasions, such as finite graphs, local finite graphs and infinite graphs. If we assume the graph is finite, then the properties of Laplacians are simple and good. But for some problems, the assumption of finite graph is obviously too narrow, so local finite graph on infinite graphs can be a better research object. We still can get good enough properties on it. In recent years, some research topics are as follows on Laplacians on infinite graph:\\
\indent (a)Definition of the operators and essential selfadjointness.\\
\indent (b)Absence of essential spectrum.\\
\indent (c)Stochastic incompleteness.\\
\indent Metric space  has a relationship with the manifold which can be seen in {\cite{6} and \cite{11}}: if we admit that there exist singular points in space, then metric space can be seen as a natural extension of the manifold. Meanwhile, it also has a similar geometric structure as the manifold. Obviously, the graph can also be seen as a kind of metric space. We can define the metric between two vertices of the graph as the natural metric, which is the number of the minimum edges connecting them. Then, we should consider whether the theories of Riemannian manifold can be extended to the graph, especially those about Ricci curvature. Many results in geometry analysis come from the Ricci curvature, especially the lower bound of Ricci curvature, such as heat kernel estimation, Harnack inequalities and Soblev inequalities. These conclusions have been made in {\cite{7}}.\\
\indent There exists an equation like this:
$$\frac{1}{2}\triangle|\bigtriangledown f|^2=\langle\bigtriangledown f,\bigtriangledown\triangle f\rangle+\|Hessf\|_2^2+Ric(\bigtriangledown f,\bigtriangledown f).$$
\indent It is an identical equation on the n-dimension Riemannian manifold. When its Ricci curvature has a lower bound, we can make a conclusion that for any $\eta\in TM$ there exists a $K\in\mathbb{R}$ that satisfies $Ric(\eta,\eta)\geqslant K|\eta|^2$.  Unfortunately in the discrete situation we can not define $\|Hessf\|_2$. But we can make use of Cauchy-Schwarz inequalities to get this inequality $\|Hessf\|_2^2\geqslant \frac{1}{n}(\bigtriangledown f)^2$. Then the Bochner inequality can be rewrite into:
$$\frac{1}{2}\triangle|\bigtriangledown f|^2\geq\langle\bigtriangledown f,\bigtriangledown\triangle f\rangle+\frac{1}{n}(\triangle f)^2+K|\bigtriangledown f|^2.$$
\indent The inequality above is the Curvature-Dimension inequality on Riemannian manifold, and we call it CD inequality for short. According to it, we can easily find that if the lower bound of the curvature is already known in the space, the "Ricci curvature" in the discrete situation can be defined. Bakery and Emery have already proved that if the chain rule is satisfied, the CD inequalities can be extended to the Markov operators on a general metric space. Yet obviously the chain rule doesn' true usually for discrete functions. Fortunately   when $p=\frac{1}{2}$, $u^p$ satisfies the chain rule even on the discrete condition. So,  {\cite{7}} introduced an improved CD inequality- CDE inequality. This definitely is a key for the research of the discrete geometry analysis.\\
\indent This paper gives an introduction of the CD inequality and several equivalent conditions of the CD inequality for the unbounded Laplacians on graph. \\
\indent The paper is organized into four parts:\\
\indent Chapter $1$ is the introduction of the graph, the Laplacians and CD inequalities on it.\\
\indent Chapter $2$ introduces some basic conclusions in order to get the main result. These conclusions include some definitions such as local finite graph, weighted graph and the domain of the operators.\\
\indent Chapter $3$ is the main conclusion of this thesis which includes some equivalent conditions of the CD inequalities.\\

\section{GRAPHS, LAPLACIANS AND CD INEQUALITIES}

\indent Given a graph $G=(V,E)$, for an $x\in V$, if there exists another $y\in V$ that satisfies $(x,y)\in E$, we call them are neighbors, and written as $x\sim y$. If there exists an $x\in V$ satisfying $(x,x)\in E$, we call it a self-loop. In this paper we allow graphs have self-loops.\\
\indent Now we will introduce some basic definitions and theorems before we get the main results.
\begin{defi} {\rm (locally finite graph)}\ \ We call a graph G is a locally finite graph if for any $x\in V$, it satisfies $\#\{y\in V|y\sim x\}<\infty$. Moreover, it is called connected if there exists a sequence $\{x_i\}_{i=0}^n$ satisfying: $x=x_0\sim x_1\sim \cdots \sim x_n=y$.\end{defi}
\begin{defi} {\rm (weighted graph)}\ \ Given a graph $G=(V,E)$, and two mappings $\mu :E\rightarrow [0,+\infty)$ and $m:V\rightarrow [0,+\infty)$ on it. $m$ is symmetric on $V$. For convenience, we extend $\mu$ onto $E$, that is to say, for any $x,y\in V$, if $x\nsim y$ or $(x,y)\not\in E$, $\mu(x,y)=0$.\end{defi}
\begin{defi} {\rm ($l^p(V,m)$ space)}\ \ Let $m$ be a measure defined as above. Then $(V,m)$ is a measure space. We define $l^p(V,m),0<p<+\infty$ space as follows:$$\{u:V \rightarrow \mathbb{R}:\sum_{x\in V}m(x)|u(x)|^p<\infty \}$$
\indent Obviously, $l^2(V,m)$ is a Hilbert space, the inner product is naturally defined as: $<\cdot,\cdot>=<\cdot,\cdot>_m$. That is: $<u,v>:=\sum_{x\in V}m(x)u(x)v(x).$ And the norm on it is defined as: $\|u\|:=<u,u>^{\frac{1}{2}}$.\end{defi}
\indent In addition, we use $l^{\infty}(V)$ to define a set including all the bounded functions on V, and we can easily know that this space is not influenced by the measure $m$. The norm on it is defined as: $\|u\|:=\sup\limits_{x\in V}|u(x)|.$
\begin{defi} {\rm (finitely supported function)}\ \ For a graph $G=(V,E)$, we call $C_0(V)$ the set of finitely supported functions if it is defined as: $C_0(V):=\{f:V\rightarrow \mathbb{R}|\# \{x\in V|f(x)\neq 0\}<\infty\}$.\end{defi}
\indent Let D is a dense subspace of $l^2(V,m)$. We define a symmetric nonnegative double mapping $Q$ on $D\times D$ to $\mathbb{R}$. D is called the domain of Q, and it is written as $D(Q)$.\\
\indent In fact, this mapping is determined by its values on the diagonal line. Then if we want to define such a mapping Q, we can just define the values on the diagonal line like this:
\begin{equation*}
Q(u):=
\left\{
\begin{aligned}
&Q(u,u)\quad &:u\in D,\\
&\infty\quad &:u\not\in D.\\
\end{aligned}
\right.
\end{equation*}
\indent If Q is lower semicontinuous, we call it closed. If Q has a closed extension it is called closable and the smallest extension is called the closure of Q as defined in {\cite{3}}.
\begin{defi} {\rm (Dirichlet form)}\ \ Q is called a Dirichlet form if it is closed and for all the contractions $C$ and $u\in l^2(V,m)$, it satisfies $Q(Cu)\leqslant Q(u)$.\end{defi}
\indent The more detailed definition can be seen in {\cite{4}}.\\
\indent On the graph the Dirichlet form has a special form as follows:$$f\mapsto Q(f):=\frac{1}{2}\sum_{x,y\in V}\mu_{xy}(f(y)-f(x))^2.$$
\indent Then we will introduce some kinds of operators on graphs.
\begin{defi} {\rm (Laplacians on locally finite graphs)}\ \ On a locally finite graph $G=(V,E,\mu,m)$ the Laplacian has a form as follows: $$\triangle f(x)=\frac{1}{m(x)}\sum_{y\in V}\mu_{xy}(f(y)-f(x)),\quad \forall f\in C_0(V).$$\end{defi}
\begin{defi} {\rm (gradient operator $\Gamma$)}\ \ The operator $\Gamma$ is defined as follows:$$\Gamma (f,g)(x)=\frac{1}{2}(\triangle(fg)-f\triangle g-g\triangle f)(x)$$.\end{defi}
\indent Always we write $\Gamma (f,f)$ as $\Gamma(f)$.
\begin{defi} {\rm (gradient operator $\Gamma_2$)}\ \ The operator $\Gamma_2$ is defined as follows:$$\Gamma_2(f,g)=\frac{1}{2}(\triangle \Gamma (f,g)-\Gamma (f,\triangle g)-\Gamma (g,\Gamma f))$$.\end{defi}
\indent Also we have $\Gamma_2(f)=\Gamma_2(f,f)=\frac{1}{2}\triangle\Gamma(f)-\Gamma(f,\triangle f)$.
\begin{defi} {\rm (nondegenerate measure)}\ \ A measure m is called nondegenerate if it satisfies $\delta :=\inf_{x\in V}m(x)>0.$\end{defi}
\indent Now we can introduce some results we need to get our main conclusions.
\begin{lem} For any $f\in l^p(V,m)$,$p\in [1,\infty)$, we have $P_tf\in l^p(V,m)$ and $$\|P_tf\|_{l^p}\leqslant \|f\|_{l^p}.$$
\end{lem}
\indent And for any $f\in l^2(V,m)$, we have $P_tf\in D(\triangle)$.
\begin{lem} For any $f\in D(\triangle)$ we have $\triangle P_tf=P_t\triangle f$.
\end{lem}
\begin{thm} Let m be a nondegenerate measure on V. Then for any $f\in l^p(V,m)$,$p\in [1,\infty)$,$$|f(x)|\leqslant \delta^{-\frac{1}{p}}\|f\|_{l^p}\quad,\forall x\in V.$$
\indent Moreover, for any $p<q\leqslant \infty$, $l^p(V,m)\hookrightarrow l^q(V,m).$
\end{thm}
\indent The proof of Lemma $2.10$, Lemma $2.11$, and Theorem $2.12$ are given by Bobo Hua and Yong Lin in {\cite{1}}.\\
\indent Now we introduce the definition of the completeness of the graph.
\begin{defi} {\rm (complete graph)}\ \ A weighted graph $(V,E,\mu,m)$ is called complete if there is a nondecreasing sequence of finitely supported functions $\{\eta\}^\infty_{k=1}$ such that $$\lim\limits_{k\rightarrow \infty}\eta_k=1\text{,and }\Gamma(\eta_k)\leqslant \frac{1}{k}.$$\end{defi}
\indent Next we will introduce two important lemmas as follows.
\begin{lem} Let $(V,E,m,\mu)$ be a complete weighted graph. Then for any $f\in D(Q)$ and $g\in D(\triangle)$,$$\sum_{x\in V}f(x)\triangle g(x)m(x)=-\sum_{x\in V}\Gamma(f,g)(x)m(x).$$
\end{lem}
\begin{lem}\
  Let $(V,E,m,\mu)$ be a complete graph. Then for any $f\in C_0(V)$ and $T>0,$ we have $\max_{[0,T]}\Gamma(P_t f)\in\ell^1_m$ and
  \begin{equation*}
  \left\|\max_{[0,T]}\Gamma(P_t f)\right\|_{\ell^1_m}\leq C_1(T,f),
  \end{equation*}
  where $C_1(T,f)$ is a constant depending on $T$ and $f$. Moveover, $$\max_{[0,T]}|\Gamma(P_t f,\frac{d}{dt} P_t f)|\in\ell^1_m\quad\mathrm{and}$$
  \begin{equation*}
    \left\|\max_{[0,T]}|\Gamma(P_t f,\frac{d}{dt} P_t f)|\right\|_{\ell^1_m}=\left\|\max_{[0,T]}|\Gamma(P_t f,\Delta P_t f)|\right\|_{\ell^1_m}\leq C_2(T,f).
  \end{equation*}
\end{lem}
\indent These two lemmas are proved in {\cite{1}}.\\
\indent Now we will introduce some basic CD inequalities(see also \cite{7} and \cite{8}).
\begin{defi} {\rm ($CD(K,\infty)$ condition)}\ \ We call a graph satisfies $CD(K,\infty)$ condition if for any $x\in V$, we have $$\Gamma_2(f)(x)\geqslant K\Gamma(f)(x),\quad K\in \mathbb{R}.$$\end{defi}
\indent For finite dimensions of curvature, we have the $CD(K,n)$ condition.
\begin{defi} {\rm ($CD(K,n)$ condition)}\ \ We call a graph satisfies $CD(K,n)$ condition if for any $x\in V$, we have $$\Gamma_2(f)\geqslant\frac{1}{n}(\triangle f)^2+K\Gamma(f).\quad K\in \mathbb{R}.$$\end{defi}
\indent Moreover, we have another condition called $CDE(x,K,n)$.
\begin{defi} {\rm ($CD(K,n)$ condition)}\ \ Let $f:V\rightarrow\mathbb{R}^+$ satisfy $f(x)>0$, $\triangle f(x)<0$. We call a graph satisfies $CD(x,K,n)$ condition if for any $x\in V$, we have$$\Gamma_2(f)(x)-\Gamma\left(f,\frac{\Gamma(f)}{f}\right)(x)\geqslant\frac{1}{n}(\triangle f)(x)^2+K\Gamma(f)(x).\quad K\in \mathbb{R}.$$\end{defi}
\indent There some relations among these conditions as follows:\\
\indent $1$.If semigroup $P_t=e^{t\triangle}$ is a diffusion semigroup, $CD(K,n)$ and $CDE'(K,n)$ are consistent.\\
\indent $2$.On graphs, $CDE'(x,K,n)$ implies $CDE(x,K,n)$ and $CD(x,K,n)$.\\
\indent Then all the preparations we need have be done. Now we will introduce the mains results of this paper.

\section{MAIN RESULTS}

\indent When we look for the equivalent properties of CD inequalities, we often set a condition: $D_{\mu}:=\max\limits_{x\in V}\frac{deg(x)}{\mu(x)}<\infty$. And the equivalent properties have already been proved in {\cite{5}} for these bounded Laplace operator on graphs. For unbounded Laplace operator, the following equivalent properties under the condition of nondegenerate measure were proved in {\cite{1}} by Bobo HUA and Yong Lin.
\begin{rem} Let $G=(V,E,m,\mu)$ be a complete graph and m is nondegenerate, i.e. $inf_{x\in V}m(x)>0$. Then the following are equivalent:\\
\indent (a) G satisfies $CD(K,\infty)$.\\
\indent (b) For any finitely supported function f,$$\Gamma(P_tf)\leqslant e^{-2Kt}P_t(\Gamma(f)).$$
\indent (c) For any $f\in D(Q)$,$$\Gamma(P_tf)\leqslant e^{-2Kt}P_t(\Gamma(f)).$$
\end{rem}
\indent In this section, similarly in {\cite{1}} we will give some equivalent properties of $CD(K,\infty)$ and $CD(K,n)$.
\begin{thm} Let $G=(V,E,m,\mu)$ be a complete graph and m is nondegenerate. Then the following are equivalent:\\
\indent (a) G satisfies $CD(K,\infty)$.\\
\indent (b) For any finitely supported function f,$$P_t(f)^2-(P_tf)^2\leqslant \frac{1-e^{-2Kt}}{K}P_t(\Gamma(f)).$$
\indent (c) For any $f\in D(Q)$,$$P_t(f)^2-(P_tf)^2\leqslant \frac{1-e^{-2Kt}}{K}P_t(\Gamma(f)).$$
\end{thm}
{\bf Proof}
\indent First, for any $f,\xi\in C_0(V)$, we make this equation
$$G(s)=\sum\limits_{x\in V}(P_{t-s}f)^2(x)P_s\xi (x)m(x)$$.
\indent Taking formal derivative of $G(s)$ in $s$, we get£º
$$G'(s)=\sum_{x\in V}(-2P_{t-s}f\triangle P_{t-s}fP_s(\xi(x))m(x)+(P_{t-s}f)^2(x)\triangle P_s\xi(x)m(x))$$
\\ \indent Now we have to show that $G(s)$ is differentiable in s.
\\ \indent For the first part£º
\begin{equation*}
\left.
\begin{aligned}
&\quad 2\sum_{x\in V}|P_{t-s}f\triangle P_{t-s}f)||P_s(\xi(x))|m(x)\\
&\leqslant 2||P_s(\xi(x))||_{l^{\infty}}||\triangle P_{t-s}f||_{l^{\infty}}(\sum_{x\in V}|P_{t-s}f|m(x))
\end{aligned}
\right.
\end{equation*}
\indent For $f, \xi\in C_0(V)$, from lemma $2.10$ we can get£º
\\ \indent $$||P_s\xi(x)||_{l^{\infty}}\leqslant ||\xi||_{l^{\infty}}<\infty$$
\indent For $f\in C_0(V)$, we know $P_{t-s}f\in D(\triangle)$ and $||\triangle P_{t-s}f||_{l^{\infty}}=||P_{t-s}\triangle f||_{l^{\infty}}\leqslant ||\triangle f||_{l^{\infty}}<\infty$
\\ \indent So we have£º
\begin{equation*}
\left.
\begin{aligned}
&\quad 2\sum_{x\in V}|P_{t-s}f\triangle P_{t-s}f)||P_s(\xi(x))|m(x)\\
&\leqslant 2||\xi ||_{l^{\infty}}||\triangle f||_{l^{\infty}}||P_{t-s}f||_{l_m^1}\\
&\leqslant 2||\xi ||_{l^{\infty}}||\triangle f||_{l^{\infty}}||f||_{l_m^1}<\infty\\
\end{aligned}
\right.
\end{equation*}
\indent For the second part, notice that $f, \xi\in C_0(V)$ and $\xi (x)\in D(\triangle)$,
\begin{equation*}
\left.
\begin{aligned}
&\quad \sum_{x\in V}(P_{t-s}f)^2(x)\triangle P_s(\xi(x))m(x)\\
&\leqslant ||\triangle P_s(\xi(x))||_{l^{\infty}}||P_{t-s}f||_{l_m^2}^2\\
&\leqslant ||P_s\triangle \xi (x)||||f||_{l_m^2}^2\\
&\leqslant ||\triangle \xi||_{l^{\infty}}||f||_{l_m^2}^2<\infty
\end{aligned}
\right.
\end{equation*}
\indent Then we know that $G(s)$ can be differentiable in $s$, and£º
$$G'(s)=\sum_{x\in V}(-2P_{t-s}f(x)\triangle P_{t-s}f(x) P_s\xi(x)m(x)+(P_{t-s}f)^2(x)\triangle P_s\xi(x)m(x)).$$
\indent For $f\in C_0(V)$, from lemma $2.10$ and Theorem $2.12$ we can easily get $(P_{t-s}f)^2\in D(Q)$,
\\ \indent Then from Lemma $2.14$, we get£º
$$G'(s)=\sum_{x\in V}(-2P_{t-s}f\triangle P_{t-s}fP_s\xi m(x)-\Gamma((P_{t-s}f)^2,P_s\xi)m(x))$$
\indent Now we replace $P_s\xi$ of $h$, and $h$ satisfies $0<h\in C_0(V)$, that is to say, $h$ is a finitely supported function.
\\ \indent Then£º
\begin{equation*}
\left.
\begin{aligned}
&\quad \sum_{x\in V}(-2P_{t-s}f\triangle P_{t-s}fh(x)m(x)-\Gamma (P_{t-s}f)^2,h(x))m(x))\\
&=\sum_{x\in V}(-2P_{t-s}f\triangle P_{t-s}fh(x)m(x)+\triangle (P_{t-s}f)^2h(x)m(x))\\
&=\sum_{x\in V}2\Gamma (P_{t-s}f)h(x)m(x)
\end{aligned}
\right.
\end{equation*}
\indent Then for $0<h\in D(Q)$, let $h_k=h\eta_k$, and $\eta_k$ satisfies
$$\lim_{k\rightarrow \infty}\eta_k=1,\Gamma (\eta_k)\leqslant \frac{1}{k},k\in N.$$
\indent we can get $0<h_k\in C_0(V)$. Then let $k\rightarrow \infty$, and for any $0<h\in D(Q)$, we have
\begin{equation*}
\left.
\begin{aligned}
&\quad \sum_{x\in V}(-2P_{t-s}f\triangle P_{t-s}fh(x)m(x)-\Gamma ((P_{t-s}f)^2,h(x))m(x))\\
&=\sum_{x\in V}2\Gamma (P_{t-s}f)h(x)m(x).
\end{aligned}
\right.
\end{equation*}
\indent For $\xi \in C_0(V)$, we easily know $P_s\xi \in D(Q)$. Then let $h=P_s\xi$, we have
\begin{equation*}
G'(s)=\sum_{x\in V}2\Gamma (P_{t-s}f)P_s\xi m(x)
\end{equation*}
\indent Integrate the equation from $0$ to $t$ by both sides£º
\begin{equation*}
\left.
\begin{aligned}
\int_0^t(\sum_{x\in V}2\Gamma (P_{t-s}f)P_s\xi m(x)) &=\int_0^tG'(s)\\
&=G(t)-G(0)\\
&=\sum_{x\in V}f^2(x)P_t\xi (x)m(x)-\sum_{x\in V}(P_tf)^2\xi (x)m(x)
\end{aligned}
\right.
\end{equation*}
\indent Since $P_t$ is a self-adjoint operator on $l_m^2$, the right hand side of the equation can be changed into:
\begin{equation*}
\left.
\begin{aligned}
\int_0^t(\sum_{x\in V}2\Gamma (P_{t-s}f)P_s\xi m(x)) &=\int_0^tG'(s)\\
&=\int_0^t\sum_{x\in V}2P_s\Gamma (P_{t-s}f)\xi (x)m(x)\\
&=\sum_{x\in V}P_t(f)^2\xi (x)m(x)-\sum_{x\in V}(P_tf)^2\xi (x)m(x)
\end{aligned}
\right.
\end{equation*}
\indent For $\xi (x)\in C_0(V)$, let $\xi (x)=\delta_y(x)$ $($when $y=x$, $\delta_y(x)=1$£¬otherwise $\delta_y(x)=0$ $)$.
\\ \indent Then, the equation is changed into:
$$P_t(f)^2(y)-(P_tf)^2(y)=2\int_0^tP_s\Gamma (P_{t-s}f)(y).$$
\indent From Remark $3.1$ we now have£º
\begin{equation*}
\left.
\begin{aligned}
P_t(f^2)-(P_tf)^2 &\leqslant2\int_0^t e^{2Ks}\Gamma (P_s\circ P_{t-s}f)ds\\
&=2\int_0^t e^{2Ks}ds \cdot \Gamma (P_tf)\\
&=\frac{e^{2Kt-1}}{K}\Gamma (P_tf)
\end{aligned}
\right.
\end{equation*}
\indent As the change in the proof is equivalent, the properties of $CD(K,\infty)$ are still equivalent properties.\\
\indent Also we can get equivalent properties of $CD(K,n)$.
\begin{thm} Let $G=(V,E,m,\mu)$ be a complete graph and m is nondegenerate. Then the following are equivalent:\\
\indent (a) G satisfies $CD(K,n)$.\\
\indent (b) For any finitely supported function f,$$\Gamma(P_tf)\leqslant e^{-2Kt}P_t\Gamma(f)-\frac{2}{n}\int_0^t e^{-2Ks}P_S(\triangle P_{t-s}f)^2ds,\quad 0<s<t.$$
\indent (c) For any $f\in D(Q)$,$$\Gamma(P_tf)\leqslant e^{-2Kt}P_t\Gamma(f)-\frac{2}{n}\int_0^t e^{-2Ks}P_S(\triangle P_{t-s}f)^2ds,\quad 0<s<t.$$
\end{thm}
{\bf Proof}
First, for any $f, \xi\in C_0(V)$, we build this functional equation£º
$$G(s)=e^{-2Ks}\sum\limits_{x\in V}\Gamma(P_{t-s}f)(x)P_s\xi(x)m(x).$$
Taking formal derivative of $G(s)$, we define the function as $A$. Then\\
\newline
$A=-2Ke^{-2Ks}\sum\limits_{x\in V}\Gamma(P_{t-s}f)(x)P_s\xi(x)m(x)+e^{-2Ks}\sum\limits_{x\in V}(-2\Gamma(P_{t-s}f,\triangle P_{t-s}f)(x)\newline P_s\xi(x)m(x)
+e^{-2Ks}\sum\limits_{x\in V}\Gamma(P_{t-s}f)\triangle P_s\xi(x)m(x).$\\
\newline
\indent Now we will show that $G(s)$ is differentiable in s.\\
\indent Without loss of generality, we assume that $\epsilon<s<t-\epsilon$ for some $\epsilon>0.$\\
\indent For the first part, from Lemma $2.15$ we have
\begin{equation*}
\left.
\begin{aligned}
&\quad |-2Ke^{-2Ks}\sum\limits_{x\in V}\Gamma(P_{t-s}f)(x)P_s\xi(x)m(x)|\\
&\leqslant|C_1|\sum\limits_{x\in V}|\Gamma(P_{t-s}f)(x)||P_s\xi|m(x)\\
&\leqslant|C_1|\|P_s\xi\|_{l^{\infty}}\|\Gamma(P_{t-s}f)\|_{l_m^1}\\
&<\infty\\
\end{aligned}
\right.
\end{equation*}
\indent $C_1$ is a constant, satisfying $|-2Ke^{-2Ks}|\leqslant C_1$.\\
\indent For the second part, from Lemma $2.15$ we have
\begin{equation*}
\left.
\begin{aligned}
&\quad |e^{-2Ks}\sum\limits_{x\in V}(-2\Gamma(P_{t-s}f,\triangle P_{t-s}f)(x)P_s\xi(x)m(x)|\\
&\leqslant |C_2|\sum\limits_{x\in V}|(-2\Gamma(P_{t-s}f,\triangle P_{t-s}f)(x)||P_s\xi(x)|m(x)|\\
&\leqslant |C_2|\|P_s\xi\|_{l^{\infty}}\|\Gamma(P_{t-s}f,\triangle P_{t-s}f)\|_{l_m^1}\\
&< \infty
\end{aligned}
\right.
\end{equation*}
\indent $C_2$ is some constant, satisfying $|e^{-2Ks}|\leqslant C_2$.\\
\indent For the last part, from Lemma $2.15$ we have
\begin{equation*}
\left.
\begin{aligned}
&\quad |e^{-2Ks}\sum\limits_{x\in V}\Gamma(P_{t-s}f)\triangle P_s\xi(x)m(x)|\\
&\leqslant |C_2|\sum\limits_{x\in V}|\Gamma(P_{t-s}f)||\triangle P_s\xi(x)|m(x)\\
&=|C_2|\sum\limits_{x\in V}|\Gamma(P_{t-s}f)||P_s\triangle \xi(x)|m(x)\\
&\leqslant |C_2|\|\triangle \xi\|_{l^{\infty}}\|\Gamma(P_{t-s}f)\|_{l_m^1}\\
&< \infty\\
\end{aligned}
\right.
\end{equation*}
\indent $C_2$ is defined as above.\\
\indent Then we can know that $G(s)$ is differentiable in $s$, and\\
\newline
$G'(s)=-2Ke^{-2Ks}\sum\limits_{x\in V}\Gamma(P_{t-s}f)(x)P_s\xi(x)m(x)+e^{-2Ks}\sum\limits_{x\in V}(-2\Gamma(P_{t-s}f,\triangle P_{t-s}f)(x)\newline P_s\xi(x)m(x)
+e^{-2Ks}\sum\limits_{x\in V}\Gamma(P_{t-s}f)\triangle P_s\xi(x)m(x).$\\
\newline
\indent From Lemma $2.14$ we get\\
\newline
$G'(s)=-2Ke^{-2Ks}\sum\limits_{x\in V}\Gamma(P_{t-s}f)(x)P_s\xi(x)m(x)+e^{-2Ks}\sum\limits_{x\in V}(-2\Gamma(P_{t-s}f,\triangle P_{t-s}f)(x)\newline P_s\xi(x)m(x)
+e^{-2Ks}\sum\limits_{x\in V}\Gamma(\Gamma(P_{t-s}f),P_s\xi(x)m(x).$\\
\newline
\indent Now we need to show that for all $h\in D(Q)$, we have\\
\begin{equation*}
\left.
\begin{aligned}
&\quad -2\sum\limits_{x\in V}\Gamma(P_{t-s}f,\triangle P_{t-s}f)(x)h(x)m(x)+\sum\limits_{x\in V}\Gamma(\Gamma(P_{t-s}f),h(x))m(x)\\
&=\sum\limits_{x\in V}\Gamma_2(P_{t-s}f)h(x)m(x).
\end{aligned}
\right.
\end{equation*}
\indent Obviously, this equation holds for all the finitely supported functions.\\
\indent Now taking a series of functions $\{\eta_k\}$ in $C_0(V)$ defined as definition $2.10$. Let $h_k=h\eta_k$, obviously $h_k\in C_0(V)$, then
\begin{equation*}
\left.
\begin{aligned}
&\quad -2\sum\limits_{x\in V}\Gamma(P_{t-s}f,\triangle P_{t-s}f)(x)h_k(x)m(x)+\sum\limits_{x\in V}\Gamma(\Gamma(P_{t-s}f),h_k(x))m(x)\\
&=\sum\limits_{x\in V}\Gamma_2(P_{t-s}f)h_k(x)m(x).
\end{aligned}
\right.
\end{equation*}
\indent Let $k\rightarrow \infty$, then for all $h\in D(Q)$, we can get
\begin{equation*}
\left.
\begin{aligned}
&\quad -2\sum\limits_{x\in V}\Gamma(P_{t-s}f,\triangle P_{t-s}f)(x)h(x)m(x)+\sum\limits_{x\in V}\Gamma(\Gamma(P_{t-s}f),h(x))m(x)\\
&=\sum\limits_{x\in V}\Gamma_2(P_{t-s}f)h(x)m(x).
\end{aligned}
\right.
\end{equation*}
\indent For $\xi\in V$, $P_s\xi\in D(Q)$, then let $h=P_s\xi$, we get
\begin{equation*}
\left.
\begin{aligned}
&\quad -2\sum\limits_{x\in V}\Gamma(P_{t-s}f,\triangle P_{t-s}f)(x)P_s\xi(x)m(x)+\sum\limits_{x\in V}\Gamma(\Gamma(P_{t-s}f),P_s\xi(x))m(x)\\
&=\sum\limits_{x\in V}\Gamma_2(P_{t-s}f)P_s\xi(x)m(x).
\end{aligned}
\right.
\end{equation*}
\indent Then $G'(s)$ can be rewritten as
$$G'(s)=e^{-2Ks}\sum\limits_{x\in V}(\Gamma_2(P_{t-s}f)-K\Gamma(P_{t-s}f))P_s\xi(x)m(x)$$
\indent By use of the equivalent properties of $CD(K,n)$, we get
$$G'(s)\geqslant e^{-2Ks}\sum\limits_{x\in V}\frac{2}{n}(\triangle P_{t-s}f)^2(x)P_s\xi(x)m(x)$$
\indent Now integrate the equation from $0$ to $t$ in $s$ by both sides, we can get
\begin{equation*}
\left.
\begin{aligned}
&\quad \int_0^tG'(s)=G(t)-G(0)\\
&=e^{-2Kt}\sum\limits_{x\in V}\Gamma(f)(x)P_t\xi(x)m(x)-\sum\limits_{x\in V}\Gamma(P_tf)(x)\xi(x)m(x)\\
&\geqslant \frac{2}{n}\int_0^te^{-2Ks}\sum\limits_{x\in V}(\triangle P_{t-s}f)^2P_s\xi(x)m(x)ds
\end{aligned}
\right.
\end{equation*}
\indent Since $P_t$ is a self-adjoint operator on $l_m^2$, we can get
\begin{equation*}
\left.
\begin{aligned}
&e^{-2Kt}\sum\limits_{x\in V}P_t\Gamma(f)(x)\xi(x)m(x)-\sum\limits_{x\in V}\Gamma(P_tf)(x)\xi(x)m(x)\\
&\geqslant \frac{2}{n}\int_0^te^{-2Ks}\sum\limits_{x\in V}P_s(\triangle P_{t-s}f)^2\xi(x)m(x)ds
\end{aligned}
\right.
\end{equation*}
\indent Let $\xi(x)=\delta_y(x)$, then
\begin{equation*}
\left.
\begin{aligned}
&e^{-2Kt}P_t\Gamma(f)(y)m(y)-\Gamma(P_tf)(y)m(y)\\
&\geqslant \frac{2}{n}\int_0^te^{-2Ks}P_s(\triangle P_{t-s}f)^2m(y)ds
\end{aligned}
\right.
\end{equation*}
\indent For $m(y)>0$,
\begin{equation*}
\left.
\begin{aligned}
&e^{-2Kt}P_t\Gamma(f)-\Gamma(P_tf))\\
&\geqslant \frac{2}{n}\int_0^te^{-2Ks}P_s(\triangle P_{t-s}f)^2ds
\end{aligned}
\right.
\end{equation*}
\indent that is to say,
$$\Gamma(P_tf)\leqslant e^{-2Kt}P_t(\Gamma(f))-\frac{2}{n}\int_0^te^{-2Ks}P_s(\triangle P_{t-s}f)^2ds.$$
\indent As the change in the proof is equivalent, the properties of $CD(K,n)$ in the theorem are still equivalent properties.

\end{document}